\documentstyle[twoside]{article}
\title{On the asymptotic expansion of $\Gamma(x)$, Lagrange's inversion theorem and the Stirling coefficients}
\author{\sc R. B. Paris\\
\\
\em{University of Abertay Dundee, Dundee DD1 1HG, UK}\\
E-Mail: r.paris@abertay.ac.uk}
\begin{document}
\def\f#1#2{\mbox{${\textstyle \frac{#1}{#2}}$}}
\def\dfrac#1#2{\displaystyle{\frac{#1}{#2}}}
\def\boldal{\mbox{\boldmath $\alpha$}}
\newcommand{\bee}{\begin{equation}}
\newcommand{\ee}{\end{equation}}
\newcommand{\lam}{\lambda}
\newcommand{\ka}{\kappa}
\newcommand{\al}{\alpha}
\newcommand{\th}{\theta}
\newcommand{\om}{\omega}
\newcommand{\Om}{\Omega}
\newcommand{\fr}{\frac{1}{2}}
\newcommand{\fs}{\f{1}{2}}
\newcommand{\g}{\Gamma}
\newcommand{\br}{\biggr}
\newcommand{\bl}{\biggl}
\newcommand{\ra}{\rightarrow}
\newcommand{\mbint}{\frac{1}{2\pi i}\int_{c-\infty i}^{c+\infty i}}
\newcommand{\mbcint}{\frac{1}{2\pi i}\int_C}
\newcommand{\mboint}{\frac{1}{2\pi i}\int_{-\infty i}^{\infty i}}
\newcommand{\gtwid}{\raisebox{-.8ex}{\mbox{$\stackrel{\textstyle >}{\sim}$}}}
\newcommand{\ltwid}{\raisebox{-.8ex}{\mbox{$\stackrel{\textstyle <}{\sim}$}}}
\renewcommand{\topfraction}{0.9}
\renewcommand{\bottomfraction}{0.9}
\renewcommand{\textfraction}{0.05}
\newcommand{\mcol}{\multicolumn}
\date{}
\maketitle
%\pagestyle{myheadings}
%\markboth{\hfill  {\it }  \hfill}
%{\hfill {\it An extension of Saalsch\"utz's theorem} \hfill}
\begin{abstract}
We show how the asymptotic expansion for the gamma function $\g(x)$, similar to that obtained by Boyd [Proc. Roy. Soc. London {\bf A447} (1994) 609--630], can be obtained by using a form of Lagrange's inversion theorem with a remainder. 
A (possibly) new closed-form representation for the Stirling coefficients is given.

\vspace{0.4cm}

\vspace{0.3cm}

\end{abstract}

\vspace{0.3cm}

\begin{center}
{\bf 1. \  Introduction}
\end{center}
\setcounter{section}{1}
\setcounter{equation}{0}
\renewcommand{\theequation}{\arabic{section}.\arabic{equation}}
The gamma function $\g(x)$ has the well-known asymptotic expansion as $x\ra\infty$
\bee\label{e11}
\g(x)=\int_0^\infty e^{-\tau}\tau^{x-1}d\tau\sim\sqrt{2\pi} x^{x-\fr}e^{-x}\sum_{n=0}^\infty \frac{(-)^n \gamma_n}{x^n},
\ee
where $\gamma_n$ are the so-called Stirling coefficients, the first few being (with $\gamma_0=1$)
\[\gamma_1=-\f{1}{12},\quad\gamma_2=\f{1}{288},\quad \gamma_3=\f{139}{51840},\quad \gamma_4=-\f{571}{2488320}.\]
%,\quad\gamma_5=-\f{163879}{209018880}.\]
The above expansion holds for large complex $x$ in the sector $|\arg\,x|\leq\pi-\delta$, $\delta>0$, although in this note we shall restrict our attention throughout to positive values of $x$. The slowly varying part of $\g(x)$ (when $x$ is large) is given by
\bee\label{e12}
\g^*(x)=\frac{\g(x)}{\sqrt{2\pi}x^{x-\fr}e^{-x}}
\ee
and, from (\ref{e11}), its asymptotic expansion is
\[\g^*(x)\sim \sum_{n=0}^\infty\frac{(-)^n \gamma_n}{x^n}=1+\frac{1}{12x}+\frac{1}{288x^2}-\frac{139}{51840x^3}+\cdots \qquad (x\ra\infty).\]

Employing the reformulation of the method of steepest descents developed by Berry \& Howls \cite{BH} (for a summary, see \cite[pp.~94--99]{P}), Boyd \cite{Bo} established the result for positive integer $m$
\bee\label{e13a}
\g^*(x)=\sum_{n=0}^{m-1}\frac{(-)^n\gamma_n}{x^n}+{\tilde R}_m(x),
\ee
where
\bee\label{e13b}
{\tilde R}_m(x)=\frac{x^{-m}}{\sqrt{2\pi}}\int_0^\infty e^{-w} w^{m-\fr}\frac{1}{2\pi i}\int_{C'}\frac{\{h(z)\}^{-m+\fr}}{h(z)-w/x}dz\,dw.
\ee
The quantity $h(z)=e^z-1-z$ and $C'$ (for $m\geq 1$) is a contour that can be taken to be a pair of straight parallel lines situated on either side of the real $z$-axis. By expanding the contour $C'$ to coincide with the other saddle points of the integrand in (\ref{e11}), Boyd then obtained the elegant expression
\[{\tilde R}_m(x)=\frac{i^mx^{-m}}{2\pi i} \int_0^\infty s^{m-1}e^{-2\pi s}\left\{\frac{\g^*(is)}{1-is/x}-(-)^m\frac{\g^*(-is)}{1+is/x}\right\}ds,\]
from which he was able to derive a bound on ${\tilde R}_m(x)$ (valid for complex $x$). This bound has been recently improved in \cite{Nem} by employing more refined bounds on $\g^*(is)$.

The Stirling coefficients appearing in the expansions (\ref{e11}) and (\ref{e13a}) can be generated numerically by means of the following recurrence relation:
\[\gamma_n=(-2)^n\frac{\g(n+\fs)}{\surd\pi}\,d_{2n},\]
\[d_n=\frac{n+1}{n+2}\left\{\frac{d_{n-1}}{n}-\sum_{j=1}^{n-1}\frac{d_j d_{n-j}}{j+1}\right\}\qquad (n\geq 1),\]
where $d_0=1$ and an empty sum is interpreted as zero.
A closed-form representation involving the 3-associated Stirling number $S_3(\ell,k)$ is found in \cite{Com} as
\[\gamma_n=\sum_{j=0}^{2n}\frac{(-)^j S_3(2j+2n,j)}{2^{j+n} (j+n)!},\]
where
\[\exp\bl[u(\frac{t^3}{3!}+\frac{t^4}{4!}+\cdots)\br]=\sum_{k, \ell\geq 0} S_3(\ell,k)\frac{u^kt^\ell}{\ell!}.\]
A proof of this result is given in \cite{BH}. A different representation has been obtained recently in \cite{LPPS}
in the form
\[\gamma_n=\sum_{m=0}^{2n}\sum_{r=0}^m\frac{(\fs)_{m+n}}{r! 2^{r-m-n}}\sum_{j=0}^{m-r}\frac{(-)^{j+n} S_{2m+2n-2r-j}^{(m-r-j)}}{j! (2m+2n-2r-j)!},\]
where $S_k^{(m)}$ denotes the Stirling number of the first kind \cite[p.~824]{AS}.

In this note we obtain the expansion of $\g^*(x)$ in the form (\ref{e13a}) and (\ref{e13b}) by making use of Lagrange's inversion theorem with a remainder, so that the inversion is valid on an infinite interval. The derivation of the remainder in Lagrange's inversion theorem is given in the appendix. The approach we use also provides a (possibly) new closed-form representation for the Stirling coefficients.
\vspace{0.6cm}

\begin{center}
{\bf 2. \ The expansion for $\g^*(x)$ as $x\ra\infty$}
\end{center}
\setcounter{section}{2}
\setcounter{equation}{0}
\renewcommand{\theequation}{\arabic{section}.\arabic{equation}}
We make the change of variable $t=\log (\tau/x)$ in Euler's integral representation for $\g(x)$
in (\ref{e11}) to find
\[\g(x)=x^xe^{-x}\int_{-\infty}^\infty e^{-xh(t)}dt,\]
where
\[ h(t)=e^t-t-1.\]
The scaled gamma function defined in (\ref{e12}) then becomes
\bee\label{e21}
\g^*(x)=\frac{x^\fr}{\sqrt{2\pi}}\int_{-\infty}^\infty e^{-xh(t)}dt.
\ee

The function $h(t)$ has saddle points (where $h'(t)=0$) at $t=2\pi ki$, $k=0, \pm 1, \pm 2, \ldots$ . The saddle at $t=0$ is the active saddle and the integration path in (\ref{e21}) coincides with the paths of steepest descent from the origin.
We now make the quadratic transformation
\bee\label{e22}
h(t)=\fs u^2
\ee
with the assumption that sign($t)=$ sign($u$),
to yield
\bee\label{e22a}
\g^*(x)=\frac{x^\fr}{\sqrt{2\pi}}\int_{-\infty}^\infty e^{-\frac{1}{2}xu^2} \frac{dt}{du}\,du.
\ee

To proceed we require the inversion of (\ref{e22}) to express $t$ as a function of the new variable $u$.
Following the argument given in \cite[p.~54]{ETC}, it can be seen that the inversion $t(u)$ is a many-valued function with branch points at $u=0$ and $u=\pm 2\sqrt{\pi k}e^{\pm\pi i/4}$, $k=1, 2, \ldots$ . Since
\[\frac{dt}{du}=\frac{u}{e^t-1},\]
the only singularities of $t(u)$ are at these branch points, and so the series expansion of $t(u)$ will converge in
$|u|<2\surd\pi$.  
\vspace{0.4cm}

\noindent 2.1 {\it The derivation of the expansion for $\g^*(x)$}
\vspace{0.2cm}

We employ Lagrange's inversion theorem with a remainder given in the appendix to obtain the inversion $t(u)$ valid for $u\in [0,\infty)$. Writing (\ref{e22}) in the form
\bee\label{e23}
u=\frac{t}{\phi(t)},\qquad \phi(t)=\left(\frac{\fs t^2}{e^t-t-1}\right)^{\!1/2}
=\left(1+2\sum_{r=1}^\infty\frac{t^r}{(r+2)!}\right)^{\!-1/2},
\ee
we have from (\ref{a3})
\[t=\sum_{n=1}^{m-1}\frac{u^n}{n!} D^{n-1} \phi^n(0)-\frac{u^m}{(m-1)!} D^{m-1} \phi^m(0)+Q_m(u)\]
for positive integer $m$, where $D^k\phi(0)\equiv (d/dt)^k \phi(t)|_{t=0}$ ($k=0,1, 2, \ldots $), and
\bee\label{e23a}
Q_m(u)=\frac{u^m}{2\pi i}
\oint_C\frac{1-u\phi'(z)}{z-u\phi(z)}\,\frac{\phi^m(z)}{z^{m-1}} dz.
\ee
The contour $C$ denotes a closed path described in the positive sense surrounding the points $z=0$ and $z=t$.
Making the change of summation index $m\ra 2m$ and differentiating we find
\bee\label{e23b}
\frac{dt}{du}=\sum_{n=0}^{m-1}\frac{u^{2n}}{(2n)!} D^{2n}\phi^{2n+1}(0)+\mbox{odd\ terms\ in\ $u$}+\frac{d}{du}Q_{2m}(u),
\ee
where we have not specified the terms in the finite sum with odd parity in $u$ since they make no contribution to the integral in (\ref{e22a}).

Substitution of the expansion (\ref{e23b}) into (\ref{e22a}) then produces
\begin{eqnarray}
\g^*(x)&=&\frac{x^\fr}{\sqrt{2\pi}} \sum_{n=0}^{m-1} \frac{D^{2n}\phi^{2n+1}(0)}{(2n)!} \int_{-\infty}^\infty u^{2n}e^{-\frac{1}{2}xu^2}du+R_m(x)\nonumber\\
&=&\frac{1}{\surd\pi}\sum_{n=0}^{m-1} \frac{2^n \g(n+\fs)}{(2n)! x^n}D^{2n}\phi^{2n+1}(0)+R_m(x),\label{e28c}
\end{eqnarray}
where the remainder after $m$ terms $R_m(x)$ is given by
\bee\label{e2r}
R_m(x)=\frac{x^\frac{3}{2}}{\sqrt{2\pi}} \int_{-\infty}^\infty ue^{-\frac{1}{2}xu^2}Q_{2m}(u)\,du.
\ee
Identification of the coefficients in the finite sum in terms of the Stirling coefficients $\gamma_n$ (see (\ref{e11})) then yields
\bee\label{e28a}
\g^*(x)=\sum_{n=0}^{m-1} \frac{(-)^n \gamma_n}{x^n}+R_m(x),
\ee
where
\bee\label{e28b}
\gamma_n=\frac{(-2)^n}{\surd\pi}\,\frac{\g(n+\fs)}{(2n)!}\,D^{2n} \phi^{2n+1}(0)
=\frac{(-)^n}{2^n n!}\,D^{2n} \phi^{2n+1}(0).
\ee
\vspace{0.4cm}

\noindent 2.2 {\it An integral representation for the remainder $R_m(x)$}
\vspace{0.2cm}

Substituting the representation of $Q_m(u)$ in (\ref{e23a}) into the expression for the remainder $R_m(x)$ in (\ref{e2r}) we obtain
\begin{eqnarray*}
R_m(x)&=&\frac{x^\frac{3}{2}}{\sqrt{2\pi}} \int_0^\infty u^{2m+1}e^{-\frac{1}{2}xu^2}\frac{1}{2\pi i}\oint_C\frac{\phi^{2m}(z)}{z^{2m-1}}\left\{\frac{1-u\phi'(z)}{z-u\phi(z)}-\frac{1+u\phi'(z)}{z+u\phi(z)}\right\}dz\,du\\
&=&\frac{2x^\frac{3}{2}}{\sqrt{2\pi}} \int_0^\infty u^{2m+2}e^{-\frac{1}{2}xu^2}\frac{1}{2\pi i}\oint_C\frac{z(-\phi(z)/z)'}{1-u^2 (\phi(z)/z)^2} (\phi(z)/z)^{2m}dz\,du.
\end{eqnarray*}
Since, from (\ref{e23}), 
\[\phi(z)/z=(2h(z))^{-1/2},\qquad (\phi(z)/z)'=-h'(z)/(2^{3/2} h^{3/2}(z)),\]
we then find after some straightforward rearrangement, together with the change of variable $w=\fs xu^2$, that
\bee\label{e25}
R_m(x)=\frac{x^{-m}}{\sqrt{2\pi}}\int_0^\infty e^{-w}w^{m+\fr}\frac{1}{2\pi i}\oint_C\frac{zh'(z)}{h(z)-w/x}\{h(z)\}^{-m-\fr}dz\,du,
\ee
where the contour $C$ denotes a closed path described in the positive sense surrounding the points $z=0$ and the two zeros (one positive and one negative) of $h(z)=w/x$.
\vspace{0.2cm}

\noindent{\bf Remark 1.}\ \ \ As in (\ref{e13b}), the contour $C$ in (\ref{e25}) can be replaced by $C'$ which is a pair of parallel lines just above and below the real $z$-axis.
\vspace{0.2cm}

\noindent{\bf Remark 2.}\ \ \ Referring to (\ref{e13b}), we see that Boyd's expression for the remainder after $m$ terms is given by
\bee\label{e26}
\frac{x^{-m}}{\sqrt{2\pi}}\int_0^\infty e^{-w}w^{m-\fr}\frac{1}{2\pi i}\oint_{C'}\frac{\{h(z)\}^{-m+\fr}}{h(z)-w/x}\,dz\,du.
\ee
We have been unable to demonstrate the equivalence between this form of the remainder and that in (\ref{e25}). We believe, however, that these two expressions are equivalent, a conjecture that is supported by high-precision numerical evaluation of the double integrals using {\em Mathematica}. In the particular case $m=2$, $x=8$ for example, we found agreement between the remainder terms in (\ref{e25}) and (\ref{e26}) to more than 30dp.

\vspace{0.6cm}

\begin{center}
{\bf 3. \ A representation for the Stirling coefficients $\gamma_n$}
\end{center}
\setcounter{section}{3}
\setcounter{equation}{0}
\renewcommand{\theequation}{\arabic{section}.\arabic{equation}}
Our representation for the Stirling coefficients is given in the following theorem.
\newtheorem{theorem}{Theorem}
\begin{theorem} $\!\!\!.$
The Stirling coefficients $\gamma_n$ $(n\geq 1)$ are given by
\bee\label{e30}
\gamma_n=2^n\sum \frac{(-2)^m (\fs)_{m+n}}{\prod_{k=1}^{2n}m_k!((k+2)!)^{m_k}},
\ee
where $(a)_n=\g(a+n)/\g(a)$ is Pochhammer's symbol,
\[m=m_1+m_2+\cdots +m_{2n}\]
and the summation is taken over all nonnegative integer solutions $(m_1, \ldots ,m_{2n})$ of the partition
\bee\label{e30a}
P_{2n}=\{(m_1, m_2, \ldots , m_{2n}): \ \sum_{k=1}^{2n}km_k=2n\}.
\ee
\end{theorem}
\vspace{0.2cm}

\noindent {\em Proof.}\ \ \ From (\ref{e28b}), the Stirling coefficients $\gamma_n$ are given by
\bee\label{e31}
\gamma_n=\frac{(-2)^n}{\surd\pi}\,\frac{\g(n+\fs)}{(2n)!}\,D^{2n} \phi^{2n+1}(0),
\ee
where $\phi(t)$ is defined in (\ref{e23}). 
To evaluate the derivatives $D^{2n} \phi^{2n+1}(0)$ we make use of Fa\`a di Bruno's formula \cite[p.~823]{AS}, \cite[p.~5]{DLMF}
\bee\label{e32}
\frac{d^{n}}{dt^{n}}[f(g(t))]=n! \sum  f^{(m)}(g(t)) \prod_{k=1}^{n} \frac{1}{m_k!} \left(\frac{g^{(k)}(t)}{k!}\right)^{\!m_k},
\ee
where
\[m=m_1+m_2+\cdots +m_{n}\]
and the summation is taken over all nonnegative integer solutions $(m_1, \ldots ,m_{n})$ of the partition
\[m_1+2m_2+ \cdots +nm_{n}=n.\]

From (\ref{e23}), we set $f(u)=u^{-n-1/2}$ and $g(t)=1+2\sum_{r=1}^\infty t^r/(r+2)!$. Then a simple calculation shows that
\[f^{(k)}(1)=(-)^k \frac{\g(n+k+\fs)}{\g(n+\fs)}, \qquad g^{(k)}(0)=\frac{2}{(k+1)(k+2)}\]
for $k=1, 2, \dots $. From (\ref{e32}) we then obtain
\[D^{2n} \phi^{2n+1}(0)=\frac{(2n)!\surd\pi}{\g(n+\fs)} \sum \frac{(-2)^m (\fs)_{m+n}}{\prod_{k=1}^{2n} m_k!\,((k+2)!)^{m_k}}.\]
Substitution of these values into (\ref{e31}) then yields the result in (\ref{e30}). 
\hfill $\Box$

An alternative version of (\ref{e30}) is
\bee\label{e30b}
\gamma_n=\frac{2^n}{(2n)!}\sum \frac{(-2)^m (\fs)_{m+n} C_{\vec m}}{\prod_{k=1}^{2n}((k+1)(k+2))^{m_k}},
\ee
where the coefficients $C_{\vec m}$ are given by\footnote{In \cite[p.~831]{AS} these quantities are called $M_3=(2n;m_1, m_2,\ldots , m_{2n})'$.}
\[C_{\vec m}=\prod_{k=1}^{2n}\frac{(2n)!}{m_k!\,(k!)^{m_k}}.\]
Values of these coefficients for $n\leq 5$ are tabulated in \cite[p.~831]{AS}.
\vspace{0.6cm}

\begin{center}
{\bf 4. \ Concluding remarks}
\end{center}
\setcounter{section}{4}
\setcounter{equation}{0}
\renewcommand{\theequation}{\arabic{section}.\arabic{equation}}
In (\ref{e28c}) we have obtained the expansion of the scaled gamma function $\g^*(x)$ as a finite sum involving inverse powers of $x$ together with a remainder $R_m(x)$ using Lagrange's inversion theorem. This result is similar to that found by Boyd \cite{Bo} who employed the Berry-Howls reformulation of the treatment of Laplace-type integrals. From this we derived an expression for the Stirling coefficients $\gamma_n$ given in (\ref{e28b}) and in Theorem 1.

A superficially similar procedure (but not equivalent) has been described by Brassesco and M\'endez \cite{BM}. They started with the result\footnote{This follows from the Euler integral for $\g(x+1)$, followed by the change of variable $\tau\ra xt$ and use of the result $\g(x+1)=x\g(x)$.} (for $x>0$)
\[\g(x)=x^x\int_0^\infty e^{-xt} t^xdt\]
and made the {\it linear} transformation $t\ra1+w$, $w=ux^{-1/2}$ to obtain
\begin{eqnarray*}
\g(x)&=&x^{x}e^{-x}\int_{-1}^\infty e^{x\{\log (1+w)-w\}}dw\\
&=&x^{x-\fr}e^{-x}\int_{-\surd x}^\infty e^{-u^2/2} e^{u^2\lambda(w)} du,
\end{eqnarray*}
where
\[\lambda(z)=z^{-2}\{\log (1+z)-z+\fs z^2\}.\]
Substituting the Maclaurin expansion
\bee\label{e41}
e^{u^2\lambda(z)}=\sum_{j\geq 0} \frac{z^j}{j!} D^j e^{u^2\lambda(z)}|_{z=0},\qquad D\equiv \frac{d}{dz}
\ee
with $z$ replaced by $ux^{-1/2}$ into the above integral, they found upon reversal of the order of summation and integration
\bee\label{e42}
\g(x)\sim x^{x-\fr}e^{-x}\sum_{j\geq 0} \frac{x^{-j/2}}{j!} D^j\left(\int_{-\surd x}^\infty u^j e^{-u^2\Lambda(z)/2}du\right)_{\!\!z=0},\quad \Lambda(z)=1-2\lambda(z).
\ee
The above integral is then extended over $(-\infty, \infty)$, so that the terms with odd index $j$ vanish, to yield the representation for the Stirling coefficients
\bee\label{e43}
\gamma_n=\frac{(-1)^n}{2^n n!} D^{2n} \Lambda^{-n-\fr}(0).
\ee
This representation is equivalent to that in (\ref{e28b}).

The implication here is that the evaluation of the $\gamma_n$ by this means has resulted in the neglect of exponentially small terms produced by extending the above integral to include the interval $(-\infty, -\surd x)$.
In addition, Brassesco and M\'endez \cite[Eq. (2.26)]{BM} incorrectly write (\ref{e42}) as an equality when this cannot be the case since the expansion (\ref{e41}) is convergent in $|z|<1$. This fact results in integration
of the series on $[0,\infty)$ beyond its interval of convergence.
In our treatment, we make the quadratic transformation in (\ref{e22}) to obtain the Stirling coefficients expressed {\it exactly} in terms of an integral over the interval $(-\infty, \infty)$. This results in {\it no exponentially small terms being neglected.} Also the use of the Lagrange inversion theorem with a remainder
circumvents the problem of integration beyond the interval of convergence (which in the case of $t(u)$ in (\ref{e23}) is $|u|<2\surd\pi$) and leads to an expression for the remainder term in the expansion.
 
The closed-form expression for the Stirling coefficients $\gamma_n$ in (\ref{e30}), and its alternative form (\ref{e30b}), involves the partition $P_{2n}$. The cardinality of this set is equal to the partition function $p(n)$, where $p(n)$ represents the number of partitions of the positive integer $n$. 
To illustrate the use of (\ref{e30b}) we take the case $n=2$, so that $p(4)=5$ and \cite[p.~831]{AS}
\[P_4=\{(0,0,0,1), (1,0,1,0), (0,2,0,0), (2,1,0,0), (4,0,0,0)\},\]
\[C_{\vec m}=\{1,4,3,6,1\}.\]
Then
\[\gamma_2=\frac{2^2}{4!}\left\{-\frac{2\cdot 1(\fs)_3}{5\cdot 6}+4(\fs)_4\left(\frac{4}{2\cdot 3\cdot 4\cdot 5}+\frac{3}{(3\cdot 4)^2}\right)-\frac{8\cdot 6(\fs)_5 }{(2\cdot 3)^2 3\cdot 4}+\frac{16\cdot 1(\fs)_6}{(2\cdot 3)^4}\right\}=\frac{1}{288}.\]
It is clear that $p(n)$ grows rapidly with $n$. Consequently, (\ref{e30b}) is not a practical means for the computation of these coefficients for large values of $n$.

\vspace{0.6cm}

\begin{center}
{\bf Appendix: \ The Lagrange expansion theorem with a remainder}
\end{center}
\setcounter{section}{1}
\setcounter{equation}{0}
\renewcommand{\theequation}{\Alph{section}.\arabic{equation}}
Let $f(t)$ and $\phi(t)$ be analytic on and inside a simple closed contour $C$ in the complex $t$-plane 
surrounding the point $t=a$. Suppose further that the function $\psi(z)=z-a-u\phi(z)$ has only one root $z=t$ inside $C$ given by
\bee\label{a0}
t-a=u\phi(t),
\ee
where $u$ is the expansion variable. The procedure we adopt is a modification of that presented in \cite[p.~17]{GW}.

Our starting point is the identity
\[f(t)=\frac{1}{2\pi i}\oint_C f(z)\frac{\psi'(z)}{\psi(z)}dz=\frac{1}{2\pi i}\oint_Cf(z)\frac{1-u\phi'(z)}{z-a-u\phi(z)}dz.\]
Upon expansion of the factor $(z-a-u\phi(z))^{-1}$ as a finite geometric progression of $m$ terms with a remainder, we find
\[
f(t)=\frac{1}{2\pi i}\oint_C f(z)(1-u\phi'(z))\left\{\sum_{n=0}^{m-1}\frac{u^n\phi^n(z)}{(z-a)^{n+1}}+\frac{u^m\phi^m(z)}{(z-a)^m (z-a-u\phi(z))}\right\}dz.
\]
Making use of the Cauchy formula
\[F^{(n)}(z)=\frac{n!}{2\pi i}\oint\frac{F(\zeta)}{(\zeta-z)^{n+1}}d\zeta,\]
we obtain
\begin{eqnarray*}
f(t)&=&\sum_{n=0}^{m-1}\frac{u^n}{n!}D^n[f(a)\phi^n(a)(1-u\phi'(a))]+Q_m(u)\\
&=&\sum_{n=0}^{m-1}\frac{u^n}{n!} D^{n}[f(a)\phi^n(a)-\frac{u}{n+1}f(a)D\phi^{n+1}(a)]+Q_m(u),
\end{eqnarray*}
where $D\equiv d/da$, the remainder $Q_m(u)$ is given by
\bee\label{a1}
Q_m(u)=\frac{u^m}{2\pi i}\oint_Cf(z)\frac{1-u\phi'(z)}{z-a-u\phi(z)}\,\frac{\phi^m(z)}{(z-a)^m}dz,
\ee
and the points $z=a$ and $z=t$ are enclosed by the contour $C$.

Straightforward rearrangement of the sum over $n$ then yields Lagrange's expansion with a remainder in the 
form\footnote{We note that the usual form of this theorem \cite[p.~17]{GW}, \cite[p.~133]{WW} has the additional requirement
$|u\phi(z)|<|z-a|$
for points on $C$, so that the arbitrary function $f(t)$ then has the expansion
\[f(t)=f(a)+\sum_{n=1}^{\infty}\frac{u^n}{n!}D^{n-1}[f'(a)\phi^n(a)].\]}
\bee\label{a2}f(t)=f(a)+\sum_{n=1}^{m-1}\frac{u^n}{n!} D^{n-1}[f'(a)\phi^n(a)]-\frac{u^m}{m!} D^{m-1}[f(a)D\phi^m(a)]+Q_m(u)
\ee
for positive integer $m$, where $\phi(t)$ is specified by (\ref{a0}).

In the special case $f(t)=t$ and $a=0$, we have from (\ref{a1}) and (\ref{a2}) the expansion for positive integer $m$
\bee\label{a3}
t=\sum_{n=1}^{m-1}\frac{u^n}{n!} D^{n-1} \phi^n(0)-\frac{u^m}{(m-1)!} D^{m-1} \phi^m(0)+\frac{u^m}{2\pi i}
\oint_C\frac{1-u\phi'(z)}{z-u\phi(z)}\,\frac{\phi^m(z)}{z^{m-1}} dz,
\ee
where $\phi(t)$ is specified in (\ref{a0}) and we have used the fact that 
\[D^{m-1}[a D\phi^m(a)]_{a=0}=(m-1) D^{m-1}\phi^m(0)\] 
and the contour $C$ encloses the poles at $z=0$ and $z=t$.

\end{document}